\documentstyle[12pt]{article}
\def\E{{\rm I\kern-.1567em E}}
\def\N{{\rm I\kern-.1567em N}}
\def\R{{\rm I\kern-.1567em R}}
\def\Pro{{\rm I\kern-.1567em P}}
\def\C{{\rm C\kern-6.5pt  
\vrule height 7.7pt width 0.4pt depth -0.5pt \phantom {.}}}                      
\def\Q{{\rm Q\kern-6.5pt  
\vrule height 7.7pt width 0.4pt depth -0.5pt \phantom {.}}}
\def\L1{L_1(\mu ,Y)}

\def\Ls*{(L_{\infty}(\mu ,X^{*}),\sigma ^{*})}

\def\t*{(L_{\infty}(\mu ,X^{*}),\tau ^{*})}
\def\P1{(\pi_n)_{n\in\N}}

\def\d1{L_1^{w^*}(\mu,Y^*)}

\def\l1{(L_1(\mu ,X),\sigma ')}

\newcommand{\fkt}[5]{\begin{array}{rrcl}
                         #1: & #2 &\rightarrow&#3\\
                     \mbox{ } & #4 &\mapsto    &#5 \end{array}}
\newcommand{\Ref}[1]{(\ref{#1})}
\newcommand{\mdE}{|\;}
\newcommand{\bgl}{\begin{eqnarray}}
\newcommand{\bglst}{\begin{eqnarray*}}
\newcommand{\egl}{\end{eqnarray}}
\newcommand{\eglst}{\end{eqnarray*}}
\newcommand{\norm}[1]{\left\|#1\right\|}

\newcommand{\ball}{{\cal B}_}
\newtheorem{theo}{Theorem}[section]

\newtheorem{prop}[theo]{Proposition}
\newtheorem{rem}[theo]{Remark}

\newcommand{\Pel}{Pel\-czy\'ns\-ki}
\newcommand{\Grabstein}{\rule{1.2ex}{1.2ex}}
\newcommand{\ebew}{\makebox[1mm]{} \hfill\Grabstein}
\newcommand{\eing}[1]{|_{#1}}
\newcommand{\eps}{\varepsilon}
\newcommand{\injB}{\otimes_{\epsilon}}
\newcommand{\projB}{\otimes_{\pi}}
\begin{document}
\author{}
\begin{center}
{\large\bf Factorization of completely bounded weakly compact operators}\bigskip\\
{\sc Hermann Pfitzner and  Georg Schl\"uchtermann}
\end{center}
\bigskip

\begin{abstract}
We prove that in the setting of operator spaces the result of Davis, Figiel, Johnson and \Pel\
on factoring weakly compact operators holds accordingly.
Though not related directly to the main theorem we add a remark on the description of
weakly compact subsets in the dual of noncommutative vector valued ${\rm L}_1$. 
\end{abstract}
\section{Introduction and preliminaries}
In 1974 Davis, Figiel, Johnson and \Pel\ \cite{DFJP} proved that a weakly
compact operator
between two Banach spaces factors through a reflexive Banach space.
(See also, e.g. \cite[p. 227]{Die-Seq})
In this note we adapt this result
to the setting of operator spaces (Th.\ \ref{theofact}).
While on the Banach space level we simply repeat the well known construction
of \cite{DFJP}, we will use some results of Pisier \cite{Pi-OH,Pi-ncLp} in order to keep trace of the operator space structures.

Pisier adapts the complex interpolation method for Banach spaces in a canonical manner to
operator spaces by constructing a (canonical) operator space structure on the usual Banach
interpolation space.

On the one hand this interpolation method serves to introduce the direct 
$l_p$-sum $l_p (X_i)$ in the operator space sense of a family $(X_i)_{i\in I}$
of operator spaces $X_i \subset {\rm B}(H_i)$. (As usual, $l_p(X_i)$ stands for the
space of families $(e_i)_{i\in I}$ with $e_i\in X_i$ for all $i$ and with
$\norm{(e_i)}=(\sum \norm{e_i}^p)^{1/p}<\infty$ if $1\leq p<\infty$ and
$\norm{(e_i)}=\sup \norm{e_i}$
if $p=\infty$. For a normed space $X$ we denote by ${\rm B}(X)$ the space of linear, bounded operators from $X$ to $X$.)
On $l_\infty (X_i)$ one defines a canonical operator space structure via its embedding in
${\rm B}(\bigoplus_2 H_i)$ which can be described by ${\rm M}_{n}(l_\infty (X_i))=
l_\infty ({\rm M}_{n}(X_i))$.
On $l_1 (X_i)$ one defines an operator space structure via the embedding $l_1 (X_i)\subset
(l_\infty (X_i^*))^*$ (because $(l_1 (X_i))^*=l_\infty (X_i^*)$ as Banach spaces), where the dual of $l_\infty (X_i^*)$ bears, of course, its standard operator space structure in the sense of \cite{Ble-stand}.
Now both $l_1 (X_i)$ and $l_\infty (X_i)$ embed continuously in the topological product $\Pi X_i,$ and it is possible to consider $l_p (X_i)$ as the complex interpolation space $(l_1 (X_i),l_\infty (X_i))_\theta$ with $\theta=1/p,$ which can be endowed with
  an operator space structure by ${\rm M}_{n}(l_p (X_i))=({\rm M}_{n}(l_1 (X_i)),
{\rm M}_{n}(l_\infty (X_i))_\theta$, see \cite{Pi-OH} for details.

On the other hand, the complex interpolation method serves, in a similar manner as just described for the direct $l_p$-sums, to construct ``noncommutative vector valued Schatten spaces'' $S_p[X]$:
If $S_p$ (respectively $S_p^n$) denote the Schatten classes
of operators on the Hilbert space $H=l_2$ (respectively on $H=l_2^n$)
and if $X$ is any operator space then the ``$X$-valued $S_\infty$'' and the
``$X$-valued $S_1$'' are defined as operator spaces by
\bgl
S_\infty [X]=S_\infty \otimes_{\rm min} X,\,\,\,\,\,
   S_1[X]=S_1 \hat{\otimes} X;  \label{gl6}
\egl
here $S_\infty$ (=compact operators on $H$) has its natural operator space structure, $S_1$ is
the dual of $S_{\infty}$  and bears
the standard dual operator space structure (\cite{Ble-stand,ER-new}, and 
$\otimes_{\rm min}$  and $\hat{\otimes}$ denote the minimal tensor product and the
projective operator space tensor product (see \cite{BlePau,ER-approx,ER-new,ER-Haag}).

We cite two results of Pisiers'  on these constructions.
The operator space structure of any operator space $Y$ can be computed by
\bgl
\norm{(y_{i,j})}_{{\rm M}_{n}(Y)} = \sup\left\{\norm{a\cdot (y_{i,j}) \cdot b}_{S_{p}^{n}[Y]}\mdE
     a,b\in \ball{S_{2p}^{n}}\right\}, \label{gl2}
\egl
for all $1\leq p <\infty$,
see Lemma 1.7 of \cite{Pi-ncLp}; where ``$\cdot$'' denotes the usual matrix product and  $\ball{X}$ is the unit ball of $X$ for a normed space $X.$ 

As to the direct $l_p$-sums we have 
\bgl
S_p^n[l_p(X_i)]=l_p(S_p^n[X_i]) \label{gl3}
\egl
for $1\leq p< \infty$,
see end of \S 2 of \cite{Pi-ncLp}.

\section{Factorization theorem}
As usual, for operator spaces $X,Y$ we denote by ${\rm CB}(X,Y)$ the space of linear and
completely bounded operators from $X$ to $Y$.
An operator space is called reflexive if it is reflexive as a Banach space.
\begin{theo}  
Weakly compact completely bounded operators factor through reflexive operator spaces.

More precisely, let $X,Y$ be operator spaces, let $T\in {\rm CB}(X,Y)$ be weakly
compact. Then there exist a reflexive operator space $R$ and operators
$T_1\in {\rm CB}(X,R),\
T_2\in {\rm CB}(R,Y)$,
such that $T=T_2\circ T_1$.
\end{theo}\label{theofact}
{\it Proof:}
On the Banach space level of the proof we adopt the well known construction
of $R$ i.e. we will define a sequence of equivalent norms on $Y$, and $R$ will turn out to
be the diagonal of the direct $l_2$-sum of these renormed $Y$.

For $m,n\in \N$  we define
\bglst
K_n:=\overline{T^{(n)}\bigl(\ball{M_n(X)}\bigr)},
\;\;\;\;\;\;\;
K_{n,m}=2^m K_n + 2^{-m}\ball{{\rm M}_{n}(Y)}.
\eglst
($T^{(n)}$ stands for the map ${\rm M}_n(X)\ni (x_{i,j})\mapsto (Tx_{i,j})\in {\rm M}_n(Y)$.)
On each ${\rm M}_{n}(Y)$ we use the Minkowski functional with respect to $K_{n,m}$ in order
to define new norms
$\norm{(y_{i,j})}_{n,m}=\inf\left\{\lambda>0\mdE (y_{i,j})\in\lambda K_{n,m}\right\}$
for each $m\in \N$.
The new norms are well defined on ${\rm M}_{n}(Y),$ because for arbitrary $y_{i,j}\in Y$ we have $(y_{i,j})\in\lambda 2^{-m}\ball{{\rm M}_{n}(Y)}\subset\lambda K_{n,m}$ with
$\lambda=n^2 2^m \max_{i,j}(\norm{y_{i,j}}),$ whence $\norm{(y_{i,j})}_{n,m}<\infty$.
We identify $Y$ and ${\rm M}_{1}(Y)$ set
$Y_{m}=(Y,\norm{\cdot}_{1,m})$ and ${\rm M}_{n}(Y_m)=({\rm M}_{n}(Y),\norm{\cdot}_{n,m})$.
Then each $Y_m$ is isomorphic to $Y$ as a Banach space because
\bgl
2^{-m}\ball{{\rm M}_{n}(Y)}\subset K_{n,m}
\subset (2^m \norm{T}_{\rm cb}+1) \ball{{\rm M}_{n}(Y)} \label{gl1}
\egl
holds for each $\N$ and in particular for $n=1$.
For each $m\in\N$ the norms $\norm{\cdot}_{n,m}$ on
${\rm M}_{n}(Y_m)$ yield an operator space structure on
$Y_m$ as can be checked by Ruan's characterisation \cite{Ruan1}.
[For the sake of completeness here are the details.  Let $(y_{i,j})\in
{\rm M}_n (Y_m)$ be arbitrary and let $\lambda>0$ be such that
$(y_{i,j})\in \lambda K_{n,m}$.
Then for $a,b\in {\rm M}_n$ we have that 
\bgl
a\cdot (y_{i,j})\cdot b&\in&
 \lambda(a\cdot K_{n,m}\cdot b)
=\lambda\big(2^m T^{(n)} (a\cdot \ball{{\rm M}_{n}(X)}\cdot b) + 2^{-m}a\cdot \ball{{\rm M}_{n}(Y)}
\cdot b\big) \nonumber\\
&\subset&
\lambda\norm{a}\;\norm{b} (2^m T^{(n)}\ball{{\rm M}_{n}(X)}+ 2^{-m} \ball{{\rm M}_{n}(Y)})
\label{gl11}\\
&=&
\lambda\norm{a}\;\norm{b} K_{n,m}\nonumber
\egl
whence $\norm{a\cdot (y_{i,j})\cdot b}_{n,m}
\leq \norm{a}\;\norm{(y_{i,j})}_{n,m}\;\norm{b}$.
(For (\ref{gl11}) we used that $Y$ is an operator space.)
In order to prove
$\norm{(y_{i,j})\oplus(z_{r,s})}_{n+k,m}\\ =\max(\norm{(y_{i,j})}_{n,m},\norm{(z_{r,s})}_{k,m})$
for $(y_{i,j})\in {\rm M}_{n}(Y_m)$, $(z_{r,s}) \in {\rm M}_k(Y_m)$
we first note that ``$\leq$'' is clear because
$(y_{i,j}) \in \lambda K_{n,m}$ and $(z_{r,s}) \in \mu K_{k,m}$ imply
$(y_{i,j}) \oplus (z_{r,s}) \in \max(\lambda,\mu) K_{n+k,m}$.
For the other inequality we define $a\in{\rm M}_{n+k}$ to be the diagonal matrix with $1$ in the
first $n$ entries and $0$ in the last $k$ entries of the diagonal, i.e.
$a= {\rm id}_n\oplus 0\in {\rm M}_{n}\oplus {\rm M}_{k}$.
Since $X$ and  $Y$ are operator spaces we have that
$a\cdot K_{n+k,m}\cdot a\subset K_{n,m}\oplus 0$.
Let $(y_{i,j})\in {\rm M}_{n}(Y_m)$ and $(z_{r,s})\in {\rm M}_{k}(Y_m)$. If $\lambda>0$ is such that
$(y_{i,j}) \oplus (z_{r,s}) \in \lambda K_{n+k,m}$ then
$(y_{i,j})\oplus 0=a\cdot ((y_{i,j})\oplus (z_{r,s}))\cdot a \in \lambda a\cdot K_{n+k,m}\cdot a$
whence $(y_{i,j})\in\lambda K_{n,m}$ whence $\norm{(y_{i,j})}_{n,m}\leq\norm{(y_{i,j})\oplus (z_{r,s})}_{n+k,m}$;
similarly,
$\norm{(z_{r,s})}_{k,m}\leq\norm{(y_{i,j})\oplus (z_{r,s})}_{n+k,m}$.] Again by (\ref{gl1}), each $Y_m$ is even completely
isomorphic to $Y$.

In the sequel we adopt Pisier's construction of $l_p$-sums, as described above,
in order to get the operator space $l_2(Y_m)$. From the original proof (\cite{DFJP, Die-Seq}) it
is well known that the ``diagonal''
\bglst
R=\left\{(y_m)\in l_2(Y_m)\mdE y_m=y_n\in Y\;\forall m,n\in \N\right\}\subset l_2(Y_m)
\eglst
is a reflexive Banach space and that the operators
\bglst
\fkt{T_1}{X}{R\,\,\,\,\,\,\,,}{x}{(Tx)}
  \fkt{\,\,\,T_2}{R}{Y}{(y)}{y}
\eglst
are well defined, continuous and satisfy $T=T_2 T_1$.

It remains to show that $T_1$ and $T_2$ are completely bounded.
First we use Pisier's results mentionned above. We have that
\bgl
\norm{((y_{i,j}^{(m)}))}_{{\rm M}_{n}(l_2(Y_m))}   
&\stackrel{\Ref{gl2}}{=}&
       \sup
       \left\{\norm{a\cdot ((y_{i,j}^{(m)})) \cdot b }_{S_2^n[l_2(Y_m)]}
       \mdE a,b\in \ball{S_4^n} \right\}\nonumber\\
&\stackrel{\Ref{gl3}}{=}&
\sup\left\{\norm{a\cdot ((y_{i,j}^{(m)})) \cdot b }_{l_2(S_2^n[Y_m])} \mdE a,b\in \ball{S_4^n}\right\}\nonumber\\
&=&
\sup\left\{\left(\sum_{m}\norm{a\cdot (y_{i,j}^{(m)}) \cdot b}
^2_{S_2^n[Y_m]}\right)^{1/2} 
\mdE a,b\in \ball{S_4^n}\right\}{\rm \,}\label{gl4}\\
&\leq&
\left( \sum_{m}\sup\left\{\norm{a\cdot (y_{i,j}^{(m)}) \cdot b }^2_{S_2^n[Y_m]} \mdE a,b\in \ball{S_4^n}\right\}
   \right)^{1/2}
\nonumber\\
&\stackrel{\Ref{gl2}}{=}&
\left( \sum_m \norm{(y_{i,j}^{(m)})}^2_{{\rm M}_{n}(Y_m)} \right)^{1/2}
\label{gl5}\\
&=&
\norm{((y_{i,j}^{(m)}))}_{l_2({\rm M}_{n}(Y_m))}. \nonumber
\egl
Let $(x_{i,j})\in\ball{{\rm M}_n(X)}$ and $y_{i,j}=Tx_{i,j}$. Then
$(y_{i,j})\in K_{n}
\subset 2^{-m} K_{n,m}=2^{-m}\ball{{\rm M}_{n}(Y_m)}$
whence
\bgl
\norm{(y_{i,j})}_{{\rm M}_n(Y_m)}\leq 2^{-m}   \label{gl6a}
\egl
for all $m\in \N$.
As to our notation, note that for
$((y_{i,j}^{(m)}))=((y_{i,j}^{(m)})_{m\in\N})_{i,j\leq n} \in {\rm M}_{n}(R)$
we have $y_{i,j}^{(m)}=y_{i,j}\in Y$ for all $m\in\N$. It follows that
\bglst
\norm{T_1^{(n)}((x_{i,j}))}    =
\norm{((y_{i,j}))}_{{\rm M}_n(R)}
\stackrel{(\ref{gl5})}{\leq}
\sum_m \norm{(y_{i,j})}^2 _{{\rm M}_n(Y_m)} \stackrel{(\ref{gl6a})}{\leq} 1,
\eglst
whence $\norm{T_1}_{\rm cb}\leq 1$.

Now we turn to $T_2$. Here we have 
\bglst
\norm{((y_{i,j}))}_{{\rm M}_{n}(R)}
&\stackrel{\Ref{gl4}}{=}&
\sup\left\{\left(\sum_{m}\norm{a\cdot (y_{i,j}) \cdot b}
^2_{S_2^n[Y_m]}\right)^{1/2} 
\mdE a,b\in \ball{S_4^n}\right\}                   \\
&\geq&
\sup\left\{\norm{a\cdot (y_{i,j}) \cdot b}_{S_2^n[Y_1]} 
\mdE a,b\in \ball{S_4^n}\right\}                 \\
&\stackrel{\Ref{gl2}}{=}&
\norm{(y_{i,j})}_{{\rm M}_{n}(Y_1)}                 \\
&\stackrel{\Ref{gl1}}{\geq}&
(2\norm{T}_{\rm cb}+1)^{-1} \norm{(y_{i,j})}_{{\rm M}_{n}(Y)}
\eglst
whence $\norm{T_2}\leq (2\norm{T}_{\rm cb}+1)$.
This ends the proof.\ebew\noindent\bigskip\smallskip\\
{\bf Remark:} Theorem \ref{theofact} remains true if we exchange ``weakly compact'' by ``Asplund'' (respectively by ``conditionally weakly compact'') and the reflexive space $R$
by an Asplund space (respectively by a space not containing a copy of $l_1$).
In the proof one only has to use $c_0$-sums instead of $\ell_2$-sums (cf. \cite[Th.~5.3.7]{Bou} respectively
\cite[p.~237]{Die-Seq}).
\bigskip\bigskip\noindent\\
\section{A remark on weak compactness in the dual of non-commutative vector-valued $L_1$-spaces}
In \cite{EW} a collection $K=(K_n)$ of nonvoid sets $K_n\subset {\rm M}_n(X)$, $X$ a
linear vector space, is called an absolutely matrix convex set on $X$ if it satisfies the following two conditions for all $n,m,r\in\N$:
\bgl
\alpha^* K_r \beta  &\subset &   
K_n {\rm\;\;for\;\, all \;\;} \alpha,\beta\in {\rm M}_{r,n},\
\|\alpha\|,\ \|\beta\|\le 1\label{gl10aa}\\
K_n \oplus K_m &  \subset &  K_{n+m}.\label{gl10ba}
\egl
For the following remarks let $K=(K_n)$ be a matrix convex set on a linear
vector space $X$.
\begin{rem}
\begin{itemize}
\item[(a)] If $X$ is an operator space then $(\ball{M_n(X)})$ is an example of an absolutely
matrix convex set
\item[(b)] Let $X,Y$ be linear vector spaces, let $T:X\rightarrow Y$ be linear. The image $L=(T^{(n)}(K_n))$ of $K$ under $T$ is a matrix convex set on $Y$.
\end{itemize}
\end{rem}
Let ${\cal W}$ be a von Neumann algebra, ${\cal W}_*$ its predual and let
$X$ be an operator space.
Analoguously to (\ref{gl6}) Pisier \cite{Pi-ncLp}
defines the continuous noncommutative $X$-valued
${\rm L}_1$ by
$$
{\rm L}_1({\cal W},X):={\cal W}_*\hat{\otimes} X.
$$
Then we know from \cite[Prop.~5.4]{BlePau} that
$$
L_1({\cal W},X)^*={\rm CB}(X,{\cal W})={\rm CB}({\cal W}_*,X^*).
$$
Note that for a measure space $(\Omega,\Sigma,\mu)$ and a Banach space $X$ the dual of the
Bochner space ${\rm L}_1(\mu,X)$ can be identified with the space ${\rm L}_{\infty}(\mu,X^*,X)$ of equivalence classes of $w^*$-measurable and essentially bounded functions $f:\Omega\rightarrow X^*$, which in turn is isometrically isomorphic to ${\rm B}(X,{\rm L}_{\infty}(\mu))$. Thus ${\rm CB}(X,{\cal W})$
seems to be a natural candidate for the noncommutative counterpart of
${\rm L}_{\infty}(\mu,X^*,X)$.\smallskip

(Off the category of operator spaces one may also define noncommutative vector valued
${\rm L}_1$-spaces only within the category of Banach spaces by the Banach projective tensor
product $\projB$:
For a Banach space $X$ and a von Neumann algebra ${\cal W}$ one defines the $X$-valued
noncommutative ${\rm L}_1$-space by
$
{\rm L}_1 ({\cal W},X)_{\rm{Ban}}={\cal W}_* \projB X.
$
It then seems natural to define the noncommutative vector valued ${\rm L}_{\infty}$-space by
${\rm L}_{\infty}({\cal W},X)_{\rm Ban}={\cal W}\injB X$ where $\injB$ is the Banach injective
tensor product and to define the corresponding ${\rm L}_p$-spaces by interpolation;
as to our knowledge this has not been treated in the literature. In this setting the dual of ${\rm L}_1({\cal W},X)_{\rm Ban}$ is isometrically isomorphic to ${\rm B}(X,{\cal W})={\rm B}({\cal W}_*,X^*)$ and seems therefore
to be a natural counterpart of ${\rm L}_{\infty}(\mu,X^*,X)$, too.)\smallskip

Up to now the characterization of weakly compact subsets of
${\rm L}_{\infty}(\mu,X^*,X)$ or of ${\rm L}_{\infty}(\mu,X)$, $X$ a Banach space, has not been achieved in a final satisfactory way; here we would like to generalize what so far has been obtained in \cite[Th.~2.4]{Schl-weak}
for weakly compact sets in
${\rm L}_{\infty}(\mu,X^*,X)$ to weakly compact sets in ${\rm CB}(X,{\cal W})$.
But since we do not dispose of a counterpart of $\Omega$ in the noncommutative setting, we must imitate the ideas of \cite[Th. 2.4.]{Schl-weak} vaguely by taking the set of pure states of ${\cal W}$ instead of $\Omega$.  
The following result could be obtained by applying 2.1 and  using
the method of \cite{Schl-weak} analogously. In order to exhibit an alternative way we shall not use the factorization theorem, in contrast to \cite{Schl-weak}.

For $v=(v_{i,j})\in {\rm M}_p(V)$, $w=(w_{k,l})\in {\rm M}_q(W)$, $V,W$ vector spaces, we use the notation of \cite{ER-new}, where $v\otimes w$ is an element of ${\rm M}_{pq}(V\otimes W),$ which can be described as follows: $v\otimes w$ is a $p\times p$-matrix whose $(i,j)$-th entry is a $q\times q$-matrix whose $(k,l)$-th entry is $v_{i,j}\otimes w_{k,l}$.
Given a pairing $\langle\cdot,\cdot\rangle:V\times W\rightarrow \C$
we also use the matrix pairing
$\langle\cdot,\cdot\rangle:{\rm M}_p(V)\times {\rm M}_q(W)\rightarrow {\rm M}_{pq}$
defined by $\langle v,w\rangle =(\langle v_{i,j}, w_{k,l}\rangle)_{i,j\leq p;
k,l\leq q}$.
For the notation in \Ref{gl12} below note that ${\rm M}_n(c_0)=({\rm M}_n\oplus {\rm M}_n\oplus\cdots)_{c_0}$
completely isometrically.
\begin{prop} 
Let ${\cal W}$ be  a von Neumann algebra and let $E_n$ denote the extreme points of the set of completely positive
contractions in ${\rm CB}({\cal W},{\rm M}_n)={\rm M}_n({\cal W}^*)$. (Thus $E_1$ is the set of pure states of ${\cal W}$.)\\
a)
Let $X$ be an operator space, and let
$(f_n)\subset {\rm L}_1({\cal W},X)^*={\rm CB}({\cal W}_*,X^*)$ be bounded. Then
$f_n\to 0$ weakly if and only if
there exists a weakly compact absolutely matrix convex set $K=(K_m)$ on $c_0$ such that
\bgl
(\langle f_n^{(p)}(\omega),x\rangle)_{n\in\N}\in K_{pq}
\label{gl12}
\egl
for all $\omega=(\omega_{i,j})\in E_p$, $x=(x_{k,l})\in {\rm M}_q(\ball{X})$.\\
b)
Let $X$ be a Banach space and let $(f_n)\subset {\rm L}_1 ({\cal W},X)_{\rm Ban}^*={\rm B}({\cal W}_*,X^*)$ be bounded.
Then $f_n\rightarrow 0$ weakly if and only if there exists a weakly compact absolutely convex set
$K\subset c_0$ such that
\bgl
(\langle f_n(\omega),x\rangle)_{n\in\N}\in K \label{gl12b}
\egl
for all pure states $\omega$ of ${\cal W}$ and all $x\in \ball{X}$.
\end{prop}\label{theochar}
{\it Proof:} (a) We define a linear bounded operator $T:l_1\rightarrow {\rm L}_1({\cal W},X)^*$
by $Te_n=f_n$ where $(e_n)$ is the canonical basis of $l^1$.
By elementary Banach space theory we know that the set $\{f_n\}$ is relatively weakly compact if and only if $T$ is weakly compact and we know that $f_n\rightarrow 0$ weakly if and only if $T$ is weakly compact and the range of $T^*\eing{{\rm L}_1({\cal W},X)}$ lies in $c_0$.

For any $\phi\in\ball{{\cal W}^*}$ and $x\in \ball{X}$ we consider
$\phi\otimes x$
as a linear functional on ${\rm CB}({\cal W}_*,X^*)$ with 
$\norm{\phi\otimes x}\leq 1$.
Identifying ${\rm L}_1({\cal W},X)^*$ with ${\rm CB}({\cal W}_*,X^*)$ it makes sense to 
write 
\bgl
\langle T^*(\phi\otimes x),e_n\rangle =
\langle \phi\otimes x,Te_n\rangle = 
\langle f_n(\phi),x\rangle,             \label{gl7}
\egl
i.e. $T^*(\phi\otimes x)=(\langle f_n(\phi),x\rangle)_{n\in\N}$.

Suppose now that $f_n\rightarrow 0$ weakly. Then $T^*$ is weakly compact and $c_0$-valued.
Let $K_m$ be the norm closure of
$T^{*(m)}(\ball{{\rm M}_m({\rm L}_1({\cal W},X)^{**})})$. Then $K=(K_m)\subset {\rm M}_m(c_0)$ is weakly compact
and absolutely matrix convex because $(\ball{{\rm M}_m({\rm L}_1({\cal W},X)^{**})})$ is absolutely  matrix convex.
Let $\omega=(\omega_{i,j})\in E_p$, $x=(x_{k,l})\in\ball{{\rm M}_q(X)}$.
Then $\norm{\omega\otimes x}_{{\rm L}_1 ({\cal W},X)^{**}}\leq 1$\\ and
by (\ref{gl7})
\bglst
(\langle f_n^{(p)}(\omega),x\rangle)_{n\in\N}
=T^{*\,(pq)}(\omega\otimes x)
=(T^*(\omega_{i,j}\otimes x_{k,l}))_{i,j\leq p;k,l\leq q}
\in K_{pq}.
\eglst

Conversely, suppose there is a weakly compact absolutely matrix convex set
$K=(K_m)\subset {\rm M}_m(X^*)$ such that (\ref{gl12}) holds.

In order to show that $f_n\rightarrow 0$ weakly it is enough to show that
$T^*(u)\in 4\,K_1$ for each $u\in {\cal W}_*\otimes X$ with
$\norm{u}_{{\rm L}_1({\cal W},X)}< 1$ because then the map
$T^*\eing{{\rm L}_1({\cal W},X)}$ takes its values in $c_0$ and is weakly compact and so is its adjoint map $(T^*\eing{{\rm L}_1({\cal W},X)})^*=T$.

By definition \cite[(3.1)]{ER-new}
$\norm{u}=\inf\{\norm{\alpha}\norm{\phi}\norm{x}\norm{\beta}\}$ where the infimum is taken over all decompositions $u=\alpha\cdot (\phi\otimes x)\cdot\beta$ with
$\alpha\in {\rm M}_{1,pq}$, $\phi=(\phi_{i,j})\in {\rm M}_p({\cal W}_*)$,
$x=(x_{k,l})\in{\rm M}_q(X)$, $\beta\in{\rm M}_{pq,1}$.
Thus for $\norm{u}< 1$ we may choose $\alpha, \phi, x, \beta$ such that all have norm $< 1$.
By Wittstock's theorem a normalized element in ${\rm CB}({\cal W},{\rm M}_p)={\rm M}_p({\cal W}^*)$ can be written as the linear combination of four completely positive contractions. Since each $K_m$ is $w^*$-closed in
${\rm M}_m (c_0)^{**}={\rm M}_m(l_{\infty})$ it follows from
(\ref{gl12}) and (\ref{gl7}) that
$T^{*(pq)}(\ball{{\rm M}_p({\cal W}_*)}\otimes \ball{{\rm M}_q(X)})$ is contained in
$(K_{pq}-K_{pq})+\imath(K_{pq}-K_{pq})$. The latter set is contained in
$4K_{pq}$ because each $K_m$ is absolutely convex. Then
\bglst
T^*(u)=\alpha\cdot T^{*(pq)}((\phi_{i,j}\otimes x_{k,l}))_{i,j\leq p; k,l\leq q}\cdot \beta
\in \alpha\cdot 4K_{pq} \cdot\beta\subset 4K_1
\eglst
since $(K_m)$ is absolutely matrix convex on $c_0.$\\
(b) The proof works almost like the one of part (a). We define
$T:\ell_1\rightarrow {\rm L}_1 ({\cal W},X)_{\rm Ban}^*$ by $e_n\mapsto f_n$.
Then (\ref{gl7}) holds accordingly.
Now if $f_n\rightarrow 0$ weakly, $T^*$ is weakly compact and we let $K$ be the norm closure of
$T^*(\ball{{\rm L}_1 ({\cal W},X)_{\rm Ban}})$.

For the converse implication suppose a weakly compact absolutely convex set $K\subset c_0$ satisfies (\ref{gl12b}). Take $u\in {\cal W}_*\otimes X$,
$\norm{u}_{{\rm L}_1 ({\cal W},X)_{\rm Ban}}\leq 1$. Then for each $\eps >0$ there are
$\lambda_i\geq 0$, $\phi_i\in\ball{{\cal W}_*}$, $x_i\in\ball{X}$, $i=1,\ldots,k$, such that 
$u=\sum \lambda_i\,\phi_i\otimes x_i$ and $\sum\lambda_i=1+\eps$. By (\ref{gl12b}) and (\ref{gl7}) we conclude
$T^*(\frac{u}{1+\eps})= T^*(\sum\frac{\lambda_i}{1+\eps}\,\phi_i\otimes x_i)
\in (K-K)+\imath(K-K)\subset 4\,{\rm aco}(K)$
since by the theorem of Krein Smulian the absolutely convex hull of a weakly compact set is again weakly compact. This ends the proof.
\ebew

 \bigskip

\vspace{2.0\baselineskip}

\parbox[t]{20em}{Dept. de Math\'ematiques\\
UFR Sciences\\
BP 6759\\
45067 Orl\'eans Cedex 2\\
France}
\hfill
\parbox[t]{11em}{Mathematisches Institut\\der Universit\"at M\"unchen\\Theresienstr.~39\\
  D\,-\,80333 M\"unchen\\Germany}
\end{document}